\documentclass[reqno, 12pt]{amsart}
\usepackage{amsmath}

\long\def\symbolfootnote[#1]#2{\begingroup%
\def\thefootnote{\fnsymbol{footnote}}\footnote[#1]{#2}\endgroup}

\newcommand\e{\epsilon}
\newcommand\R{{\mathbb R}}

\newtheorem{teo}{Theorem}[section]

\newtheorem{remark}[teo]{Remark}

\newtheorem{coro}[teo]{Corollary}

\newtheorem{theorem}[teo]{Theorem}
\newtheorem{lemma}[teo]{Lemma}

\catcode`\@=11

\@addtoreset{equation}{section}
\@addtoreset{teo}{section}
\usepackage{amsbsy}

\newlength{\defbaselineskip}
\newcommand{\setlinespacing}[1]%
           {\setlength{\baselineskip}{#1 \defbaselineskip}}

\begin{document}
\title{Asymptotic behaviour of a semilinear elliptic system with a large exponent}
\thanks{This research was supported by
        FONDECYT 3040059}

\author{I.A. GUERRA,}
\address{ Departamento de Matematica y C. C.,
Universidad  de Santiago,
 Ca\-si\-lla 307, Co\-rreo 2, Santiago, Chile}
 \email{\tt iguerra@usach.cl}



\begin{abstract}
Consider the problem
\begin{eqnarray*}
-\Delta u &=& v^{\frac 2{N-2}},\quad  v>0\quad \mbox{in}\quad \Omega,\\
-\Delta v &=& u^{p},\:\:\:\quad  u>0\quad \mbox{in}\quad \Omega, \\
u&=&v\:\:=\:\:0 \quad \mbox{on}\quad \partial \Omega,
\end{eqnarray*}
where $\Omega$ is a bounded convex domain in $\R^N,$ $N>2,$ with
smooth boundary $\partial \Omega.$ We study the asymptotic
behaviour of the least energy solutions of this system as $p\to
\infty.$ We show that the solution remain bounded for $p$ large
and have one or two peaks away form the boundary. When one peak
occurs we characterize its location.

\end{abstract}

\maketitle

\section{Introduction}
In this article we consider the problem
\begin{eqnarray}\label{pp}
\begin{cases}
-\Delta(-\Delta u)^{(N-2)/2} = u^{p},\quad  u>0\quad \mbox{in}\quad \Omega, \\
\qquad\qquad \qquad u=\Delta u=0 \quad \mbox{on}\quad
\partial \Omega,
\end{cases}
\end{eqnarray}
where $\Omega$ is a bounded convex domain in $\R^N,$ $N>2,$ with
smooth boundary $\partial \Omega.$ We consider the so-called least
energy solutions of \eqref{pp}, obtained by the minimization
problem
$$
c_p:=\inf\limits_{v\in W^{2,\frac{N}{2}}(\Omega)\cap
W_0^{2,\frac{N}{2}}(\Omega)}\left\{(\int\limits_\Omega|\Delta
v|^{N/2}dx)^{2/N}\colon \|v\|_{p+1}=1\,\right\}
$$
By standard argument $c_p$ is achieved by a positive function
$\underline{u}_p$, which is a positive scalar multiple of a
function solving  \eqref{pp}. Let us denote such least energy
solution by $u_p.$

Problem \eqref{pp} is the particular case $q=2/(N-2)$ for the
system
\begin{eqnarray*}
-\Delta u &=& v^{q},\quad  v>0\quad \mbox{in}\quad \Omega,\\
-\Delta v &=& u^{p},\quad  u>0\quad \mbox{in}\quad \Omega, \\
u&=&v\:\:=\:\:0 \quad \mbox{on}\quad \partial \Omega,
\end{eqnarray*}
For this system the condition on $(p,q),$ given by
\begin{equation}\label{sobhyperbola}
\frac{N}{p+1}+\frac{N}{q+1}-(N-2)=0,
\end{equation}
so called {\it critical hyperbola}, describes a  the borderline
between existence and non-existence of positive solutions. In this
article, we fix $q=2/(N-2),$ that is we {\it stand} in an
asymptote, and we prove that as we increase $p$ the least energy
solution develops a peak behavior. The case $N=4$ these type of
results were shown in \cite{BEG,T1,T2} and in the case $N=2,$ the
problem reduces to one equation, and we observe a similar
behaviour, see \cite{RW1,RW2}.

More precisely, our aim is to prove the following results.

\begin{theorem}\label{main1}
Let $u_p$ the least energy solution of \eqref{pp} There exists
$C_1,C_2$ independent of $p$ such that
$$
0<C_1<\|u_p\|_{L^{\infty}(\Omega)}<C_2<+\infty
$$
for $p$ large enough.
\end{theorem}

For the next result we define
$$
w_p:=\frac{u_p}{(\int\limits_\Omega u_p^p\,dx)^\frac 2{N-2}}
$$
For a sequence $w_{p_n}$ of $w_p,$ we define the {\it blow up set}
$S$ of $\{w_{p_n},\}$ as
\begin{align*}
S:=&\{x\in \overline{\Omega}\colon \exists\:\: \mbox{a
subsequence}\quad
w_{p_n}, \\
 & \exists \{x_n\} \subset \Omega\:\: \mbox{such that}\:\: x_n\to
x\:\: \mbox{and}\:\: w_{p_n}(x_n)\to \infty\,\}.
\end{align*}
We define a {\it peak point} $P$ for $u_p$ to be a point in
$\overline{\Omega}$ such that $u_p$ does not vanish in the
$L^\infty$ norm in any neighborhood of $P$ as $p\to \infty.$ We
shall see later that peaks point of $\{u_p\}$ are contained  in
the blow up set $S$ of $\{w_p\}$

\begin{theorem}\label{main2}
Let $\Omega$ a convex bounded domain in $\R^N,$ $N>3,$ with smooth
boundary $\partial \Omega$. Then for any sequence $w_{p_n}$ of
$w_p,$ with $p_n\to \infty$ there exists a sequence still denoted
by $w_{p_n}$ such that the blow up set $S$ of this subsequence is
contained in $\Omega$ and has the property $1\leq
\textrm{card}(S)\leq 2.$

If $\textrm{card}(S)=1$ and $S=\{x_0\}$ then: \\
1)
$$
f_n:=\frac{u_{p_n}^{p_n}}{\int\limits_\Omega
u_{p_n}^{p_n}\,dx}=(\int\limits_\Omega u_{p_n}^{p_n}\,dx)^{\frac
2{N-2}p_n-1} w_{p_n}^{p_n}\to \delta_{x_0}.
$$
in the sense of distributions.

2) $w_{p_n}\to \tilde G(\cdot,x_0)$ in
$C_{loc}^2(\overline{\Omega}\setminus\{x_0\})$ where $\tilde
G(x,y)$ solves
$$
-\Delta \tilde G(x,\cdot) = G^{\frac 2{N-2}}(x,\cdot)\quad\mbox{
in}\quad \Omega,\qquad\tilde
G(x,\cdot)=0\quad\mbox{on}\quad\partial \Omega,$$ where $-\Delta
G(\cdot,y) = \delta_{y}$ in $\Omega,$ $G(\cdot,y)=0$ on $\partial
\Omega. $

3) $x_0$ is a critical point of 
$\tilde \phi(x):=\tilde g(x,x)$ where the function $\tilde g(x,y)$
is given by
$$
\tilde g(x,y)= \tilde G(x,y)+\frac 1{(N-2)^{\frac
N{N-2}}\omega_{N-1}^{\frac 2{N-2}}}\log|x-y|,
$$
where $\omega_{N-1}$ the area of the unit sphere $S^{N-1}$ in
$\R^N.$
\end{theorem}
We observe that regularity of  $\tilde \phi$ is needed to compute
its critical points in 3).
Indeed, by definition of $\tilde G$, we have
\begin{eqnarray}\label{limtilde}
\lim\limits_ {y\to x}|x-y|^{4-N}\Delta\tilde
g(x,y)=-\frac{2}{N-2}\frac{g(x,x)}{((N-2)\omega_{N-1})^{\frac{4-N}{N-2}}}
\end{eqnarray}
for $x\in \Omega,$ where $g(x,y)$ is the regular part of $G(x,y)$,
i.e
$$
g(x,y)=G(x,y)-\frac 1{(N-2)\omega_{N-1}|x-y|^{N-2}}.
$$
By elliptic regularity, for $N\geq 3,$ the function $g(x,\cdot)$
is regular and so is $\tilde \phi.$
\begin{remark}
We conjecture that for $N=3$ the conclusions in Theorem
\ref{main2} also hold. The only difficulty is to prove Lemma
\ref{cotawn} for $N=3,$ but we think that is only technical.
\end{remark}

\section{Estimates for $c_p$}

\begin{lemma}\label{ldt}
For every $t\geq N/2.$ there is $D_t$ such that
$$
\|u\|_{t}\leq D_t t^{\frac{N-2}{N}}\|\Delta u\|_{N/2}
$$
where
$$
\lim\limits_{t\to\infty}D_t=(\frac{N-2}{Neb_0})^{\frac{N-2}{N}}
$$
with $b_0=\frac
N{\omega_{N-1}}[4\pi^{N/2}/\Gamma((N-2)/2)]^\frac{N}{N-2}=N(N-2)^{\frac
N{N-2}}\omega_{N-1}^{\frac 2{N-2}}.$
\end{lemma}

\begin{proof} From \cite{A}, we have the following Higher-Order
version of the Moser-Trudinger inequality,
$$
\int\limits_{\Omega}\exp(b_0|u|^{N/(N-2)})dx \leq C|\Omega|
$$
for any $u$ such that $\|\Delta u\|_{N/2}\leq 1.$  Therefore
\begin{align}
\frac 1{\Gamma(\frac{t(N-2)}{N}+1)}\int\limits_\Omega
u^t\,dx&=\frac 1{\Gamma(\frac{t(N-2)}{N}+1)}\int [b_0(\frac
u{\|\Delta
u\|_{N/2}})^{\frac N{N-2}}]^{\frac{t(N-2)}{N}}\,dx \\
&\times b_0^{-\frac{t(N-2)}{N}}\|\Delta u\|^t_{N/2}
\\
&\leq \int\limits_\Omega \exp(b_0(\frac u{\|\Delta
u\|_{N/2}})^{\frac N{N-2}})\,dx \times
b_0^{-\frac{t(N-2)}{N}}\|\Delta u\|_{N/2}^t
\end{align}
hence
$$
\|u\|_{L^t(\Omega)}\leq
\Gamma(t(N-2)/N+1)^{1/t}C^{1/t}b_0^{-(N-2)/N}|\Omega|^{1/t}\|\Delta
u\|_{L^{N/2}(\Omega)}
$$
We conclude using the Stirling's formula,
$$
\Gamma(\frac{t(N-2)}{N}+1)^{1/t}\sim
(\frac{N-2}{Ne})^{\frac{N-2}{N}}t^{\frac{N-2}{N}}.
$$

\end{proof}

\begin{lemma}\label{lcp}
\begin{equation}
\lim\limits_{p\to \infty}c_p p^{\frac{N-2}{N}}=(\frac{N b_0
e}{N-2})^{\frac{N-2}{N}}
\end{equation}
\end{lemma}
\begin{proof} Let $L$ such that $B_L\subset \Omega,$ and
$l\in(0,L)$ to be fixed later. Let  $m_l(x)=H((\log L/l)^{-1}\log
1/|x|),$ a regularized version of a Moser's function, where $H$ is
such that for $\e\in (0,1/2)$
$$
H(t)=\begin{cases}\e \Phi(t/\e) & 0<t\leq \e\cr t, & \e<t\leq 1-\e
\cr 1-\e \Phi((1-t)/\e) & 1-\e<t\leq 1, \cr 1 & 1<t
\end{cases}
$$
with $\Phi\in C^\infty[0,1],$  $\Phi(0)=\Phi'(0)=0$
$\Phi(1)=\Phi'(1)=1.$ Clearly $m_l\in W_0^{2,N/2}(\Omega)$ and
$m_l(x)=1$ for $|x|\in l.$  A calculation gives
$$
|\Delta m_l(x)|=\left|-\frac{a_0}{\omega_{N-1}}H'((\log
L/l)^{-1}\log 1/|x|)(\log L/l)^{-1}|x|^{-2}+O(\log
L/l)^{-2}|x|^{-2})\right|
$$
where $a_0=(\omega_{N-1}b_0/N)^{\frac{N-2}{N}}.$ Thus
$$
\int\limits_{B}|\Delta
m_l|^{N/2}dx=M^\frac{N}{2}:=\omega_{N-1}^{1-\frac{N}{2}}a_0^{\frac
N{2}}(\log 1/r)^{1-\frac N{2}}A,
$$
where $A\leq 1+C\e,$ see \cite{A} for details. We define $
\psi_l=m_l/M$ and find
$$
\left(\int\limits_\Omega \psi_l^{p+1}\,dx\right)^{\frac
1{p+1}}\geq \left(\int\limits_{B_l} \psi_l^{p+1}\,dx\right)^{\frac
1{p+1}}=\frac 1{M}(\frac 1{N}l^N\omega_{N-1})^{\frac 1{p+1}}.
$$
Take $l=L\exp(-(N-2)(p+1)/N^2)$ and recall that
$M^{-1}=(\omega_{N-1}\log L/l)^{\frac{N-2}{N}}A^{-\frac
2{N}}a_0^{-1},$ we find
$$
\int\limits_\Omega \psi_l^{p+1}dx)^{\frac 1{p+1}}\geq
\frac{\omega_{N-1}^{\frac{N-2}{N}}}{a_0}(\frac{N-2}{N^2})^{\frac{N-2}{N}}(p+1)^{\frac{N-2}{N}}A^{-\frac
2{N}}(\frac 1{N}L^N\omega_{N-1})^{\frac 1{p+1}}e^{-\frac{N-2}{N}}.
$$
Then
\begin{align} c_p(p+1)^{\frac{N-2}{N}} &\leq
a_0\omega_{N-1}^{-\frac{N-2}{N}}(\frac{N-2}{N^2})^{-\frac{N-2}{N}}A^{\frac
2{N}}(\frac
1{N}L^N\omega_{N-1})^{-\frac 1{p+1}}e^{\frac{N-2}{N}}\\
& \leq b_0^{\frac{N-2}{N}}(\frac{N
e}{N-2})^{\frac{N-2}{N}}A^{\frac 2{N}}(\frac
1{N}L^N\omega_{N-1})^{-\frac 1{p+1}}
\end{align}
Letting $p\to \infty $ and $\e\to 0,$ we obtain the result by
combining this with Lemma \ref{ldt}.
\end{proof}
\begin{coro}\label{coroconv}
We have
$$
p^{\frac{N-2}{2}}\int\limits_\Omega|\Delta
u_p|^{N/2}dx=(\frac{Nb_0e}{N-2})^{\frac{N-2}{2}},\quad
p^{\frac{N-2}{2}}\int\limits_\Omega
u_p^{p+1}dx=(\frac{Nb_0e}{N-2})^{\frac{N-2}{2}}.
$$
\end{coro}

\begin{proof}[Proof of Theorem \ref{main1}]
Let $\Lambda_N$ be given by the minimization problem
$$
\Lambda_N=\inf\left\{\frac{\|\Delta u\|_{L^\frac
N{2}(\Omega)}}{\|u\|_{L^{\frac N{2}}(\Omega)}}\mid u\in W^{2,\frac
N{2}}(\Omega)\cap W_0^{1,\frac N{2}}(\Omega)\,\right\}.
$$
From Lemma \ref{ldt}, we find  $0<\Lambda_N<\infty.$ Using this,
$$
\int\limits_{\Omega}u_p^{p+1}\,dx=\int\limits_\Omega|\Delta
u_p|^{\frac N{2}}\,dx\geq
\Lambda_N^{N/2}\int\limits_{\Omega}u_p^{N/2}\,dx.
$$
Thus $\int\limits_{\Omega}(u_p^{p+1}-\Lambda_N^{N/2}u^{N/2}_p)\geq
0$, therefore $\|u_p\|^{p+1-\frac N{2}}_{L^{\infty}(\Omega)}\geq
\Lambda_N^{\frac N{2}}.$ If $p>(N-2)/2$ then
$$
\|u_p\|_{L^{\infty}(\Omega)}>\Lambda_N^{\frac N{2(p+1-\frac
N{2})}}\geq C_1>0.
$$
To obtain an upper bound for $\|u_p\|_{L^{\infty}(\Omega)},$ let
\begin{equation}
\gamma_p=\max\limits_{x\in \Omega},\quad  \Omega_l:=\{x\in
\Omega\colon t< u_p(x)\}, \quad \mathcal{A} =\{x\in \Omega\colon
\frac{\gamma_p}{2}<u_p(x)\}.
\end{equation}
Applying Lemma \ref{ldt} and \ref{coroconv}, we obtain
$$
(\int\limits_\Omega
u_p^{\frac{N^2p}{2(N-2)}}dx)^{\frac{2(N-2)}{N^2p}}\leq
D_{\frac{N^2p}{2(N-2)}}({\frac{N^2p}{2(N-2)}})^{\frac{N-2}{N}}\|\Delta
u_p\|_{L^{\frac N{2}}(\Omega)}\leq M
$$
where $M$ is independent of $p$ for $p$ large. This implies
\begin{equation}\label{cota1}
(\frac{\gamma_p}{2})^{\frac{N^2p}{2(N-2)}}|\mathcal{ A}|\leq
M^{\frac{N^2p}{2(N-2)}}.
\end{equation}
Taking  $v_p^{\frac 2{N-2}}=-\Delta u_p$ and integrating by parts
$$
\int\limits_{\partial\Omega_t}|\nabla u_p|\,ds =
\int\limits_{\Omega_t}v_p^{\frac 2{N-2}}dx.
$$
By Coarea formula we have
$$
-\frac d{dt}|\Omega_t|=\int\limits_{\partial\Omega_t}\frac
1{|\nabla u_p|}\,ds
$$
Then Schwartz inequality gives
$$
-\frac d{dt}|\Omega_t|\int\limits_{\Omega_t}v_p^{\frac
2{N-2}}\,dx=\int\limits_{\partial\Omega_t}|\nabla
u_p|\,ds\int\limits_{\partial\Omega_t}\frac 1{|\nabla
u_p|}\,ds\geq |\partial \Omega_t|^2.
$$
The isoperimetric inequality in $\R^N$
$$
|\partial \Omega_t|\geq C_N|\Omega_t|^{\frac{N-1}{N}},
$$
yields
$$
-\frac d{dt}|\Omega_t|\int\limits_{\Omega_t}v_p^{\frac
2{N-2}}\,dx\geq C_N^2|\Omega_t|^{\frac{2(N-1)}{N}}.
$$
We define $r(t)$ such that $|\Omega_t|=\omega_{N-1}r^N(t)/N.$ Then
$$
\frac d{dt}|\Omega_t|=\omega_{N-1}r^{N-1}(t)r'(t)<0.
$$
Thus
\begin{align*}
\frac{N}{\omega_{N-1}^{\frac{N-2}{N}}C_N^2r(t)^{N-1}}&\int\limits_{\Omega_t}v_p^{\frac
2{N-2}}\,dx\geq -\frac 1{r'(t)} \\
-\frac{dt}{dr}\leq
\frac{N}{\omega_{N-1}^{\frac{N-2}{N}}C_N^2r(t)^{N-1}}\int\limits_{\Omega_t}v_p^{\frac
2{N-2}}\,dx &\leq
\frac{N}{\omega_{N-1}^{\frac{N-2}{N}}C_N^2r(t)^{N-1}}(\sup\limits_{\Omega}v_p)^{\frac
2{N-2}}|\Omega_t| \\
&=(\sup\limits_{\Omega}v_p)^{\frac 2{N-2}}
\frac{r(t)\omega_{N-1}^{\frac{2}{N}}}{NC_N^2}.
\end{align*}
Integrating this inequality from $r=0$ to $r=r_0,$ we have
$$
t(0)-t(r_0)\leq
\frac{\omega_{N-1}^{\frac{2}{N}}}{2NC_N^2}r_0^2(\sup\limits_{\Omega}v_p)^{\frac
2{N-2}}
$$
Choosing $r_0$ such that $t(r_0)=\gamma_p/2$ that is
$|\mathcal{A}|=|\Omega_{\gamma_p/2}|=\omega_{N-1}r_0^N,$ the last
inequality yields
$$
\gamma_p\leq \frac{\omega_{N-1}^{\frac{2}{N}}}{NC_N^2}r_0^2
(\sup\limits_{\Omega}v_p)^{\frac 2{N-2}}
$$
But $v_p^{\frac 2{N-2}}=-\Delta u_p$ satisfies
$$
-\Delta v_p=u_p^p\quad \mbox{in}\quad \Omega, \quad v_p=0\quad
\mbox{on}\quad \partial\Omega.
$$
By elliptic regularity (\cite{GT}, Theorem 3.7),
$\sup\limits_{\Omega}v_p\leq C\sup\limits_{\Omega}u_p^p\leq
C\gamma_p^p$ where $C=C(\Omega).$ Thus we have
\begin{equation}\label{cota2}
\gamma_p\leq \bar C
\omega_{N-1}^{\frac{2}{N}}\gamma_p^{\frac{2p}{N-2}}r_0^2=\bar
C\gamma_p^{\frac{2p}{N-2}}|\mathcal{A}|^{\frac{2}{N}}
\end{equation}
where $\bar C=C^{\frac{2p}{N-2}}/(NC_N^2).$ By \eqref{cota1} and
\eqref{cota2},
$$
\gamma_p\leq \bar C\gamma_p^{-p}(2M)^{\frac{Np}{N-2}}
$$
which implies
$$
\gamma_p\leq \bar C^{\frac 1{p+1}}(2M)^{\frac{Np}{(N-2)(p+1)}}.
$$
Therefore there exists $C>0$ independent of $p$ such that $
\gamma_p\leq  C$ for $p$ large.
\end{proof}

Next we have the corollary.

\begin{coro}
There exist $C_1,C_2>0$ independent of $p$ such that
$$
\frac{C_1}{p^{\frac{N-2}{2}}}\leq
\int\limits_{\Omega}u_p^p\,dx\leq \frac{C_2}{p^{\frac{N-2}{2}}}
$$
for large $p.$
\end{coro}
\begin{proof} From Corollary \ref{coroconv} and Theorem \ref{main1}, we have for $p$ large
$$
C'\leq p^{\frac{N-2}{2}}\int\limits_{\Omega}u_p^{p+1}\,dx\leq
\|u_p\|_{L^\infty(\Omega)}p^{\frac{N-2}{2}}\int\limits_{\Omega}u_p^p\,dx\leq
C''p^{\frac{N-2}{2}}\int\limits_{\Omega}u_p^p\,dx,
$$
where $C',C''>0$ constant independent of $p.$ This shows the left
inequality. Now by Holder inequality,
$$
p^{\frac{N-2}{2}}\int\limits_{\Omega}u_p^p\,dx\leq
(p^{\frac{N-2}{2}}\int\limits_{\Omega}u_p^{p+1}\,dx)^{\frac
p{p+1}}p^{\frac{N-2}{2(p+1)}}|\Omega|^{\frac 1{p+1}}.
$$
Using Corollary \ref{coroconv}, for $p$ large the RHS is bounded.
\end{proof}

\section{Proof of Theorem \ref{main2}}
Let
$$
w_p:=\frac{u_p}{(\int\limits_\Omega u_p^p\,dx)^{\frac
2{N-2}}}=\frac{u_p}{\lambda_p},\quad
\lambda_p:=(\int\limits_\Omega u_p^p\,dx)^{\frac 2{N-2}},
$$
and
$$
f_p(x):=\frac{u_p^p}{\int\limits_\Omega u_p^p\,dx}.
$$
This yields
\begin{eqnarray*}
-\Delta(-\Delta w_p)^{(N-2)/2} &=& f_p\quad  w_p>0\quad \mbox{in}\quad \Omega, \\
w_p&=&\Delta w_p\:\:=\:\:0 \quad \mbox{on}\quad \partial \Omega,
\end{eqnarray*}

Since $\Omega$ is convex, we can derive standard uniform boundary
estimates of $\{w_p\}$  which leads to conclude that the blow up
set of $\{w_p\}$ is contained in the interior of $\Omega.$

Using the methods in \cite{QS} Proposition 3.2, we can show that
$$
a)\:\: \int\limits_\Omega w_p\phi_1\leq C \quad \mbox{and}\quad b)
\:\:\int\limits_\Omega (-\Delta w_p)\phi_1\leq C,
$$
where $\phi_1$ is the positive eigenvalue of
$(-\Delta,H_0^1(\Omega)),$ normalized to $\int\limits_{\Omega}
\phi_1=1.$ Combining inequality $a)$ with the ideas of \cite{DLN}
based on the method of the moving planes from \cite{GNN}, we
obtain a uniform bound in the boundary. Indeed, we can find
$\delta>0$ such that for any $x\in \Omega_\delta:=\{z\in \bar
\Omega\colon d(z,\partial\Omega)<\delta\,\},$ we have
$$
w_p(x)\leq C(\Omega_\delta).
$$

Now we extend a known results from \cite{BM}.

\begin{lemma}\label{lbm}
Let $u$ be a regular solution of
\begin{eqnarray}
-\Delta(-\Delta u)^{\frac{N-2}{2}}=f(x)\quad \mbox{in}\quad
\Omega\subset \R^N
\\
u=\Delta u=0 \quad \mbox{on}\quad \partial \Omega,
\end{eqnarray}
where $f\in L^1(\Omega),$ $f\geq 0.$ For any $\e\in (0, b_0)$ we
have
$$
\int\limits_\Omega
\exp\left(\frac{(b_0-\delta)|u(x)|}{\|f\|_{L^1(\Omega)}^{\frac
2{N-2}}}\right)\,dx\leq \frac{b_0}{\delta}|\Omega|.
$$
\end{lemma}
\begin{proof}
We proof this by the symmetrization method. Consider the
symmetrized problem
\begin{align}
-\Delta(-\Delta U)^{\frac{N-2}{2}}=F(x)\quad \mbox{in}\quad
\Omega^* \\
U=\Delta U=0\quad \mbox{on}\quad \partial\Omega^*.
\end{align}
Here $\Omega^*$ is a ball centered at the origin with the same
volume at $\Omega,$ say $\Omega^*=B(0,R),$ and $F$ is the
symmetric decreasing rearrangement of $f$. By \cite{Ta1} and
\cite{Ta2}, we have
$$
u^*\leq U
$$
where $u^*$ is the symmetric rearrangement of $u$. Clearly $U$
satisfies
\begin{align}
-(r^{N-1}U'(r))'= r^{N-1}V^{\frac 2{N-2}}, \quad r\in (0,R) \\
-(r^{N-1}V'(r))'= r^{N-1}F(r), \quad r\in (0,R) \\
U'(0)=V'(0)=U(R)=V(R)=0.
\end{align}
Multiple integrations give,
$$
-U'(r)\leq \frac 1{(N-2)^{\frac N{N-2}}\omega_{N-1}^{\frac
2{N-2}}}\frac 1{r}\|F\|_{L^1(\Omega^*)}^{\frac 2{N-2}}.
$$
Hence
\begin{align}
|U(r)|\leq \frac 1{(N-2)^{\frac N{N-2}}\omega_{N-1}^{\frac
2{N-2}}}\|F\|_{L^1(\Omega^*)}^{\frac 2{N-2}}\log(\frac R{r}), \\
\int\limits_{\Omega^*}\exp\left[(N-\e)(N-2)^{\frac
N{N-2}}\omega_{N-1}^{\frac 2{N-2}}\frac
U{\|F\|_{L^1(\Omega^*)}^{\frac 2{N-2}}}\right]\leq
\int\limits_{B_R(0)}\exp\log(\frac R{r})^{N-\e}dr \\
= \omega_{N-1}\int\limits_0^R (\frac
R{r})^{N-\e}r^{N-1}\,dr=\e^{-1}\omega_{N-1}R^N.
\end{align}
Letting $\e(N-2)^{\frac N{N-2}}\omega_{N-1}^{\frac
2{N-2}}=\delta,$ we have
$$
\int\limits_{\Omega^*}\exp\left[(b_0-\delta)\frac
U{\|F\|_{L^1(\Omega^*)}^{\frac
2{N-2}}}\right]\leq\frac{b_0}{\delta}|\Omega|.
$$
By the properties of the symmetric decreasing functions,
$\|F\|_{L^1(\Omega^*)}=\|f\|_{L^1(\Omega)},$ and
\begin{align*}
\int\limits_{\Omega}\exp\left[(b_0-\delta)\frac{
u(x)}{\|f\|_{L^1(\Omega)}^{\frac 2{N-2}}}\right]\,dx=
\int\limits_{\Omega}\exp\left[(b_0-\delta)\frac{
u^*(x)}{\|f\|_{L^1(\Omega)}^{\frac 2{N-2}}}\right]\,dx \\
\leq \int\limits_{\Omega^*}\exp\left[(b_0-\delta)\frac
U{\|F\|_{L^1(\Omega^*)}^{\frac
2{N-2}}}\right]\leq\frac{b_0}{\delta}|\Omega|,
\end{align*}
which proves the lemma.
\end{proof}

\begin{coro}\label{cbm}
a) Let $u_n$ be a sequence of solutions of
\begin{eqnarray*}
-\Delta(-\Delta u_n)^{\frac{N-2}{2}}=V_ne^u_n\quad \mbox{in}\quad
\Omega\subset \R^N
\\
u_n=\Delta u_n=0 \quad \mbox{on}\quad \partial \Omega,
\end{eqnarray*}
such that $\|V_n\|_{L^p(\Omega)}\leq C_1,$ for some $p\in
(1,\infty),$  $\|V_n\|_{L^p(\Omega)}\leq C_1,$ and
$$
\|V_ne^{u_n}\|_{L^1(\Omega)}\leq \e_0<b_op/(p-1).
$$
Then $\{u_n\}$ is uniformly bounded in $L^\infty_{loc}(\Omega).$

b) Let $u_n$ be a sequence of solutions of $ -\Delta(-\Delta
u_n)^{\frac{N-2}{2}}=V_ne^{u_n}$ in $\Omega\subset \R^N$ with
$V_n\geq 0$ and $u_n,-\Delta u_n\geq 0$ on the boundary. Assume
for some $p\in (1,\infty)$ that
\begin{align}
&(1) \quad \|V_n\|_{L^p(\Omega)}\leq C_1, \\
&(2) \quad \|u_n\|_{L^1(\Omega)}\leq C_2
\\
&(3) \quad \|V_ne^{u_n}\|_{L^1(\Omega)}\leq \e_0<b_op/(p-1)
\end{align}
Then $\{u_n\}$ is uniformly bounded in $L^\infty_{loc}(\Omega).$
\end{coro}

\begin{proof} Part a).
Fix $\delta> 0,$ so that $b_0-\delta> \e_0(p'+\delta).$ By Lemma
\ref{lbm}, we have
$$
\int\limits_\Omega\exp[(p'+\delta)|u_n|]\leq C
$$
for some $C$ independent of $n$. Therefore $e^{u_n}$ is bounded in
$L^{p'+\delta}(\Omega),$ hence $V_ne^{u_n}$ is bounded in
$L^{1+\e_0}(\Omega).$ Then by elliptic regularity, we have $u_n$
bounded in $L^{\infty}_{loc}(\Omega).$

The part b) follows similarly. With restriction we may assume that
$\Omega=B_R(x_0)$ for some $x_0.$ We consider
\begin{eqnarray*}
-\Delta(-\Delta u_{1n})^{\frac{N-2}{2}}=V_ne^{u_n}\quad
\mbox{in}\quad \Omega\subset \R^N
\\
u_{1n}=\Delta u_{1n}=0 \quad \mbox{on}\quad \partial \Omega,
\end{eqnarray*}
and $-\Delta(-\Delta u_{2n})^{\frac{N-2}{2}}=0$ in $\Omega$ with
$u_{2n}=u_n\geq 0$ and  $-\Delta u_{2n}=-\Delta u_n\geq 0.$ By
uniqueness
$$
(-\Delta u_n)^{\frac{N-2}{2}}=(-\Delta
u_{1n})^{\frac{N-2}{2}}+(-\Delta u_{2n})^{\frac{N-2}{2}}.
$$
where each term is positive. For $N=3,$ we have $-\Delta u_n\leq
-2(\Delta u_{1n}+\Delta u_{2n})$ and for $N\geq 4,$ $-\Delta
u_n\leq -(\Delta u_{1n}+\Delta u_{2n}).$ In the last case, define
$H=u_n-u_{1n}-u_{2n},$ and by the maximum principle $H\leq 0$ in
$B_{R}(x_0).$ This  gives
$$
u_n\leq u_{1n}+u_{2n}\quad \mbox{in}\quad B_{R}(x_0),
$$
and similarly for $N=3,$ we have
$$
u_n\leq 2(u_{1n}+u_{2n})\quad \mbox{in}\quad B_{R}(x_0).
$$
Note that $u_{2n}\leq u_n$ and $u_{1n}\leq u_n$ in $B_{R}(x_0).$
Now a uniform bound for $u_{1n}$ is given by part a) and we know
by mean value theorem,
$$
\|u_{2n}\|_{L^\infty(B_{R/2}(x_0))}\leq
C\|u_{2n}\|_{L^1(B_{R}(x_0))}  \leq C\|u_{n}\|_{L^1(\Omega)}\leq
C,
$$
and the last inequality follows from the assumption (2).
\end{proof}

Let $u_n,$ $w_n,$ $\lambda_n,$ and $f_n$ denote $u_{p_n},$
$w_{p_n},$ $\lambda_{p_n},$ and $f_{p_n}.$ First we note that the
blow up set $S$ of the sequence $\{w_n\}$ is not empty. In fact,
$$
\sup\limits_{x\in\Omega} w_n(x)\geq \frac{C}{\lambda_n}\to
+\infty,
$$
by Theorem \ref{main1} and using that $p_n\lambda_n\leq C$ for $C$
independent of $p_n$ large. This also shows that the set peaks of
$\{u_n\}$ is contained in the set $S$. Since
$$
f_n(x)=\frac{u_n^{p_n}}{\int\limits_\Omega u_n^{p_n}\,dx}\in
L^1(\Omega), \quad f_n\geq 0, \quad \int\limits_\Omega f_n(x)\,dx
=1,
$$
there exists a subsequence (denoted also by $\{u_n\}$) such that
there exists a positive bounded measure $\mu$ in the set of real
bounded Borel measures in $\Omega,$ satisfying $\mu(\Omega)\leq 1$
and
$$
\int\limits_\Omega f_n\phi\to \int\limits_\Omega
\phi\,d\mu\quad\mbox{for all}\quad \phi\in C_0(\Omega).
$$
We now define the quantity
$$
L_0=\frac 1{e}\limsup\limits_{p\to\infty}p\left(\int\limits_\Omega
u_p^{p}\,dx\right)^{\frac 2{N-2}}.
$$
From the proof of Theorem \ref{main1}, we obtain $1\leq L_0\leq
Nb_0/(N-2).$

For any $\delta>0$ we call a point $x_0\in \Omega$ a
$\delta-$regular point of $\{u_n\}$ if there exists $\varphi\in
C_0(\Omega),$ $0\leq \varphi \leq 1$ with $\varphi\equiv 1$ in a
neighborhood of $x_0$ such that
$$
\int\limits_\Omega \varphi d\mu \leq
\left(\frac{b_0}{L_0+2\delta}\right)^{\frac{N-2}{2}}.
$$
We also define for $\delta>0,$ $\delta-$irregular set of a
sequence $\{u_n\}$ such that
$$
\Sigma(\delta)=\{x_0\in \Omega\colon \mbox{$x_0$ is not a
$\delta-$regular point}\:\}.
$$
Note that $x_0\in \Sigma(\delta)$ implies
$$
\mu(x_0)>\left(\frac{b_0}{L_0+2\delta}\right)^{\frac{N-2}{2}}.
$$

The next results is crucial to prove Theorem \ref{main2}.

\begin{lemma}\label{cotawn} Assume $N>3.$ Let $x_0$ be a $\delta-$regular point of a sequence
$\{u_n\}$ then $\{w_n\}$ is bounded in $L^\infty(B_{R_0}(x_0))$
for some $R_0>0.$
\end{lemma}

\begin{proof} Let $x_0$ be a $\delta-$regular point. Then there
exists  $R>0$ such that
$$
\int\limits_{B_{R}(x_0)}f_n\,dx<(\frac{b_0}{L_0+\delta})^{\frac{N-2}{2}}
$$
for $n$ sufficiently large.

Let $w_{1n}$ be solution of
\begin{eqnarray*}
-\Delta(-\Delta w_{1n})^{\frac{N-2}{2}}=f_n\quad \mbox{in}\quad
B_{R}(x_0)
\\
w_{1n}=\Delta w_{1n}=0 \quad \quad\mbox{on}\quad \partial
B_{R}(x_0),
\end{eqnarray*}
and let $w_{2n}$ be solution of
\begin{eqnarray*}
-\Delta(-\Delta w_{2n})^{\frac{N-2}{2}}=0\quad \mbox{in}\quad
B_{R}(x_0)
\\
w_{2n}=w_n,\quad \Delta w_{2n}=\Delta w_{n} \quad\quad
\mbox{on}\quad
\partial B_{R}(x_0),
\end{eqnarray*}
By the maximum principle we have $-\Delta w_{1n}>0$ and $-\Delta
w_{2n}>0, $ $w_{1n}>0,$ and $w_{2n}>0$ in $B_{R}(x_0).$

Clearly by uniqueness
$$
(-\Delta w_n)^{\frac{N-2}{2}}=(-\Delta
w_{1n})^{\frac{N-2}{2}}+(-\Delta w_{2n})^{\frac{N-2}{2}}.
$$
If $N\geq 4,$ $-\Delta w_n\leq -(\Delta w_{1n}+\Delta w_{2n}),$
then we can  define $H=w_n-w_{1n}-w_{2n},$ and by the maximum
principle $H\leq 0$ in $B_{R}(x_0).$ This gives
$$
w_n\leq w_{1n}+w_{2n}\quad \mbox{in}\quad B_{R}(x_0)
$$
Note that $w_{2n}\leq w_n$ and $w_{1n}\leq w_n$ in $B_{R}(x_0).$
The solution $w_{2n}$ is uniformly bounded near $x_0,$ in fact the
the mean value theorem gives
$$
\|w_{2n}\|_{L^\infty(B_{R/2}(x_0))}\leq
C\|w_{2n}\|_{L^1(B_{R}(x_0))}  \leq C\|w_{n}\|_{L^1(\Omega)}\leq
C,
$$
and the last inequality follows from Lemma \ref{lbm}.

So we need to bound $\{w_{1n}\}.$ We first choose $t$ such that
$t':=t/(t-1)=L_0+\delta/2.$ Since $L_0>1,$ there exists $t >1.$
Then we have
$$
\int\limits_{B_{R}(x_0)}f_n\,dx<(\frac{b_0}{L_0+\delta})^{\frac{N-2}{2}}
$$
Lemma \ref{lbm} implies
\begin{equation}\label{cota3}
\int\limits_{B_{R}(x_0)}
\exp(t'|w_{1n}(x)|)\,dx=\int\limits_{B_{R}(x_0)}
\exp((L_0+\delta/2)|w_{1n}(x)|)\,dx\leq C,
\end{equation}
where $C=C(\delta)\to \infty$ as $\delta\to 0.$

By the inequality $\log x\leq x/e$ for $x>0,$ we have
\begin{align}
\log f_n=&\log \frac{u_n^{p_n}}{\lambda_n^{\frac{N-2}{2}}}=p_n
\log \frac{u_n}{\lambda_n^{\frac{N-2}{2p_n}}}\leq
p_n\frac{u_n}{e\lambda_n^{\frac{N-2}{2p_n}}} \\
& \leq
\frac{L_0+\delta/3}{\lambda_n}\frac{u_n}{\lambda_n^{\frac{N-2}{2p_n}}}=
\frac{t'-\delta/6}{\lambda_n^{\frac{N-2}{2p_n}}}\frac{u_n}{\lambda_n}\leq
t'w_n(x).
\end{align}
The second inequality follows from the definition of $L_0$ and the
last form
$\lim\limits_{n\to\infty}\lambda_n^{\frac{N-2}{2p_n}}=1.$

Thus we get the pointwise estimate $f_n(x)<\exp(t' w_n(x)),$ which
implies
\begin{equation}\label{cota4}
(f_n\exp(-w_{1n}(x))^t<C\exp(t' w_{1n}(x))
\end{equation}
 in $B_{R/2}(x_0)$ because
$w_{2n}$ is uniformly bounded in $B_{R/2}(x_0)$ and $w_n\leq
w_{1n}+w_{2n}$ in $B_{R/2}(x_0).$

Rewrite the equation for $w_{1n}$ as
\begin{eqnarray*}
-\Delta(-\Delta
w_{1n})^{\frac{N-2}{2}}=\underline{f_ne^{-w_{1n}(x)}}e^{+w_{1n}(x)}\quad
\mbox{in}\quad B_{R}(x_0)
\\
w_{1n}=\Delta w_{1n}=0 \quad \mbox{on}\quad \partial B_{R}(x_0),
\end{eqnarray*}
Clearly $w_{1n}, -\Delta w_{1n}\geq 0$ in $B_{R}(x_0)$ We now
check the assumptions of Corollary \ref{cbm}. Let
$V_n=f_ne^{-w_n(x)},$
\begin{align}
&(1) \quad \|V_n\|_{L^t(B_{R/2}(x_0))}\leq C_1, \quad\mbox{ by \eqref{cota3} and \eqref{cota4} }\\
&(2) \quad \|w_{1n}\|_{L^1(B_{R/2}(x_0))}\leq C_2 \quad\mbox{ by Lemma \ref{lbm} }\\
&(3) \quad
\|V_ne^{u_n}\|_{L^1(B_{R/2}(x_0))}=\|f_n\|_{L^1(B_{R/2}(x_0))}\leq
\e_0<Ct/(t-1)
\end{align}
Applying Corollary \ref{cbm}, we conclude that $\{w_{1n}\}$ is
uniformly bounded in $B_{R/2}(x_0).$
\end{proof}

\begin{lemma}
$S=\Sigma(\delta)$ for any $\delta>0.$
\end{lemma}
\begin{proof} $S\subset \Sigma(\delta)$ is clear from Lemma
\ref{cotawn}.
 Now suppose that $x_0\in \Sigma(\delta)$ and
$\|w_n\|_{L^\infty(B_{R_0}(x_0))}<C$ for some $C$ independent of
$n.$
 then $f_n=\lambda_n^{p_n-1}w_n^{p_n}\to 0$ uniformly on
 $B_{R_0}(x_0),$ which implies $x_0$ is a $\delta-$regular point,
 that is $x_0\not\in \Sigma(\delta).$ Thus contradiction shows
 that for every $R>0$ we have $\lim\limits_{n\to
 \infty}\|w_n\|_{L^\infty(B_{R}(x_0))}=\infty$ at least for a
 subsequence. So $x_0\in S.$
\end{proof}

This lemma implies that
$$
1\geq \mu(\Omega)\geq \left(\frac{
b_0}{L_0+2\delta}\right)^{\frac{N-2}{2}}\textrm{card}(\Sigma(\delta))=\left(\frac{b_0}{L_0+2\delta}\right)^{\frac{N-2}{2}}\textrm{card}(S).
$$
Combining this with the estimate $L_0\leq Nb_0/(N-2),$ we have
$$
\textrm{card}(S)\leq
\left(\frac{L_0+2\delta}{b_0}\right)^{\frac{N-2}{2}}\leq
\left(\frac{Nb_0/(N-2)+2\delta}{b_0}\right)^{\frac{N-2}{2}}
$$
hence, since $S$ is not empty
$$
1\leq \textrm{card}(S)\leq \left(
\frac{N}{N-2}\right)^{\frac{N-2}{2}}<e.
$$
This proves the first part of Theorem \ref{main2}.

In the following we assume $\textrm{card}(S)=1$ with $S=\{x_0\}.$
Then $w_n(x)\leq C$ on any compact set $K\subset
\overline{\Omega}\setminus\{x_0\},$ which implies $f_n\to 0$
uniformly on compacts of $\overline{\Omega}\setminus\{x_0\}.$

Take $\varphi\in C_0(\Omega).$ For given $\e>0$ we choose $r>0$
small such that as $n\to \infty,$ we have
\begin{align*}
|\int\limits_\Omega f_n\varphi\,dx-\varphi(x_0)|\leq
\int\limits_\Omega f_n|\varphi(x)-\varphi(x_0)|\,dx \\
\leq \int\limits_{B_r(x_0)}f_n|\varphi(x)-\varphi(x_0)|\,dx
+\int\limits_{\Omega\setminus
B_r(x_0)}f_n|\varphi(x)-\varphi(x_0)|\,dx\leq \e.
\end{align*}
Therefore
\begin{equation}
f_n\to\delta_{x_0}
\end{equation}
in the sense of distributions. Let $\tilde w_n^{\frac
2{N-2}}=-\Delta w_n,$ then
$$
-\Delta\tilde w_n=f_n\quad \mbox{in}\quad \Omega,\quad \tilde
w_n=0\quad \mbox{on}\quad \partial \Omega.
$$
with $f_n\to 0$ uniformly in compact subsets of
$\overline{\Omega}\setminus\{x_0\}.$ This proves 1) of Theorem
\ref{main2}.

For 2), on any compact $K\subset
\overline{\Omega}\setminus\{x_0\},$ we have $\tilde w_n$ is
bounded and $f_n\to 0$ uniformly. By elliptic regularity there
exists a subsequence of $\tilde w_n,$ still denoted by $\tilde
w_n$ that approaches a function say $G'$ in $C^{2,\alpha}(K)$
weakly in $W^{1,q}(\Omega)$ for ($1<q<2$) and strongly in
$L^1(\Omega)$ by the compact embedding
$W^{1,q}(\Omega)\hookrightarrow L^1(\Omega).$ As in \cite{RW1},
$G'=G(\cdot,x_0).$ Since  $-\Delta w_n=\tilde w_n^{\frac 2{N-2}}$
by the convergence of $\tilde w_n,$ we have that $w_n$ converges
to $G''$ in $C^{2,\alpha}(\Omega),$ and by uniqueness $G''=\tilde
G.$

To prove 3) we use a Pohozaev identity. From \cite{M,V}, for any
$y\in \R^N,$ we have for any $\Omega'\subset \R^n,$ the following
identity
\begin{eqnarray}\nonumber
\int\limits_{\Omega'} \Delta u(x-y,\nabla v)+\Delta v(x-y,\nabla u) - (N-2)(\nabla u,\nabla v)dx=\\
\int_{\partial \Omega'} \frac{\partial u }{\partial n}(x-y,\nabla
v)+\frac{\partial v }{\partial n}(x-y,\nabla u)-(\nabla u,\nabla
v)(x-y,n)\, ds. \label{mitvan}
\end{eqnarray} Let
$\Omega'=\Omega.$ For the system $-\Delta v= u^p$ and $-\Delta u=
v^{\frac 2{N-2}}$ in $\Omega,$ the identity \eqref{mitvan} takes
the form
\begin{eqnarray}
(\frac{N}{p+1} -\bar a)\int\limits_\Omega u^{p+1}\,dx +(N-2-\bar b)\int\limits_\Omega v^{\frac N{N-2}}\,dx \\
+(N-2-\bar a-\bar b)\int\limits _\Omega(\nabla u,\nabla v)\,dx
=-\int\limits_{\partial \Omega}(\nabla u,\nabla v )(x-y,n)\,ds.
\end{eqnarray}
We choose $\bar a+\bar b=N-2$, $\bar a= 0$ and so $\bar b=N-2.$
This gives for $u_p$ and $v_p=(-\Delta u_p)^{\frac{N-2}{2}},$
\begin{eqnarray}\label{epoh}
\frac{N}{p+1}\int\limits_\Omega
u_p^{p+1}\,dx=-\int\limits_{\partial \Omega}\frac{\partial
u_p}{\partial n}\frac{\partial v_p}{\partial n}(n,x-y)\,ds
\end{eqnarray}
differentiating with respect to $y$,
$$
\int\limits_{\partial \Omega}\frac{\partial u_p}{\partial
n}\frac{\partial v_p}{\partial n}n\,ds=0.
$$
Taking $v_p\to G$ and $u_p\to \tilde G$ in $C^2(\Omega),$ as $p\to
\infty,$ we get
\begin{eqnarray}\label{vec2}
\int\limits_{\partial \Omega}(\nabla \tilde G(x,x_0),\nabla
G(x,x_0))n\,ds=0
\end{eqnarray}
On the other hand we have the following result.
\begin{lemma}\label{lcpoint} For every $x_0\in \Omega$
\begin{eqnarray}\label{vec4}
\int\limits_{\partial \Omega}(\nabla \tilde G(x,x_0) ,n)(\nabla
(\Delta \tilde G(x,x_0))^{\frac{N-2}{2}},n) n\,ds=-\nabla \tilde
\phi(x_0).
\end{eqnarray}
\end{lemma}

Hence combining
\eqref{vec2}
with \eqref{vec4}, we complete the proof of part 3) and the
Theorem \ref{main2} is proven.

\begin{proof}[Proof Lemma \ref{lcpoint}.] Let $\Omega'=\Omega\setminus B_r$ with $r>0.$ For a system
$-\Delta v= 0$ and  $-\Delta u= v^{\frac 2{N-2}}$ in $\Omega',$
the identity \ref{mitvan}, takes the form
\begin{eqnarray}\nonumber
\int\limits_{\Omega'} (N-2)v^{\frac N{N-2}}- \bar a v^{\frac
N{N-2}}\, dx=
\int\limits_{\partial \Omega'} \frac{N-2}{N}v^{\frac N{N-2}}(x-y,n)\,ds\\
+\int\limits_{\partial \Omega'} \frac{\partial u }{\partial
n}\left[(x-y,\nabla v)+\bar a v \right]+\frac{\partial v
}{\partial n}\left[(x-y,\nabla u)+\bar b u\right]- (\nabla
u,\nabla v)(x-y,n)\,ds \label{Mident}
\end{eqnarray}
with $\bar a+\bar b=N-2.$ We choose $\bar a=N-2$ and take
$v=G(x,0)$ and $u=\tilde G(x,0).$ Upon differentiation with
respect to $y$, \eqref{Mident} transforms into
\begin{eqnarray*}
 \int\limits_{\partial \Omega}\frac{\partial \tilde G }{\partial n}\frac{\partial G }{\partial
 n}n\,ds
=\int\limits_{\partial B_r}\left\{ \frac{N-2}{N}G^{\frac
N{N-2}}n+\frac{\partial \tilde G }{\partial n}\nabla G\, +
\frac{\partial G }{\partial n}\nabla \tilde G-(\nabla \tilde
G,\nabla G)n\,\right\} ds.
\end{eqnarray*}
Note that $u=v=0$ on $\partial \Omega,$ implies $\nabla u= (\nabla
u,n)n$ and $\nabla v= (\nabla v,n)n$ on $\partial \Omega.$ Let
$\Gamma=\omega_{N-1}(N-2),$ we have
\begin{eqnarray*}
\nabla \tilde G =-\frac 1{\Gamma^{\frac
2{N-2}}(N-2)}|x|^{-2}x+\nabla \tilde g,\quad \nabla G =-\frac
1{\omega_{N-1}}|x|^{-N}x+\nabla g,
\end{eqnarray*}
\begin{eqnarray*}
\frac{\partial \tilde G }{\partial n}=-\frac 1{\Gamma^{\frac
2{N-2}}(N-2)}|x|^{-1}+(\nabla \tilde g,n),\quad \frac{\partial G
}{\partial n}=-\frac 1{\omega_{N-1}}|x|^{1-N}+(\nabla g,n)
\end{eqnarray*}
\begin{eqnarray*}
(\nabla \tilde G,\nabla G)=\frac{|x|^{-N}}{\omega_{N-1}
\Gamma^{\frac 2{N-2}}(N-2)}-\frac{(\nabla
g,x)}{\Gamma^p(N-2)}|x|^{-2}
-\frac {(\nabla \tilde g,x)}{\omega_{N-1}}|x|^{-N}+(\nabla \tilde
g,\nabla g)
\end{eqnarray*}
and
$$
\frac{N-2}{N}G^{\frac N{N-2}}=\frac{N-2}{N}\left[\frac
{1}{\Gamma^{\frac 2{N-2}}}|x|^{-2}-\Delta \tilde
g\right]\left[\frac 1{\Gamma}|x|^{2-N}+g\right]
$$

Using $\int\limits_{\partial B_r} n=0,$ we get
\begin{eqnarray}\nonumber
 \int\limits_{\partial \Omega}\frac{\partial \tilde G }{\partial n}\frac{\partial G }{\partial
 n}n\,ds=\frac{N-2}{Nr^{N-1}}\int\limits_{\partial B_r}\left\{\frac 1{\Gamma^p}r^{N-2-1}g-\Delta \tilde g
\frac 1{\Gamma} r-\Delta \tilde g g r^{N-1}\right\}n\,ds \\
\nonumber +\frac 1{r^{N-1}}\int\limits_{\partial B_r}\{(\nabla
\tilde g,n)\nabla g+(\nabla g,n)\nabla \tilde g-(\nabla \tilde g,\nabla g)n\}r^{N-1}\,ds \\
-\frac 1{r^{N-1}}\int\limits_{\partial
B_r}\left\{\frac{1}{\omega_{N-1}}\nabla \tilde
g+\frac{r^{N-2}}{\Gamma^p(N-2)}\nabla g\right\}\,ds. \label{Grto0}
\end{eqnarray}
Since $N\geq 3,$ and $\tilde g$ and $g$ are regular, we obtain in
the limit as $r\to 0,$
\begin{eqnarray*}
 \int\limits_{\partial \Omega}\frac{\partial \tilde G }{\partial n}\frac{\partial G }{\partial
 n}n\,ds=\lim\limits_{r\to 0}\frac 1{r^{N-1}}\int\limits_{\partial B_r}\frac{1}{\omega_{N-1}}\nabla
 \hat
g\,ds=
\nabla \tilde \phi(0).
\end{eqnarray*}
\end{proof}

\end{document}